\documentclass[12pt,thmsa]{article}
\usepackage{amsfonts}
\usepackage{sw20bams}

\input tcilatex
\begin{document}

\author{Steven Finch}
\title{A\ Convex Maximization Problem:\ Discrete Case}
\date{December 5, 1999}
\maketitle

\begin{abstract}
We study a specific convex maximization problem in \textit{n}-dimensional
space. The conjectured solution is proved to be a vertex of the polyhedral
feasible region, but only a partial proof of local maximality is known.
Integer sequences with interesting patterns arise in the analysis, owing to
the number theoretic origin of the problem.
\end{abstract}

\begin{center}
\textit{Dedicated in memory of Jill Spurr Titus, with love}
\end{center}

\section{Problem}

For each positive integer $n$, maximize the convex function 
\[
\sum_{i=1}^n\frac 1{x_i} 
\]
over the polyhedron in $n$-dimensional real space $\Bbb{R}^n$ defined by

\[
\begin{array}{cccc}
(j+1)x_j+x_i\geq (j+1)i+\varepsilon _{i\,j} &  & \text{for} & 1\leq j\leq
i\leq n
\end{array}
\]
where $\varepsilon _{i\,j}=1$ if $i=j=1$ and $\varepsilon _{i\,j}=0$
otherwise. Prove that:

\begin{enumerate}
\item[(i)]  a global maximum $\left( a_1,a_2,...,a_n\right) $ exists and is
unique

\item[(ii)]  the components $a_i$ of the global maximum satisfy 
\[
\begin{array}{ccccc}
a_1=1, &  & a_2=2, &  & a_3=4
\end{array}
\]
and, when $i\geq 4$, 
\[
a_i=(j+1)(i-a_j)
\]
\end{enumerate}

for any $j$ with $(j+1)a_j-ja_{j-1}\leq i<(j+2)a_{j+1}-(j+1)a_j$ .\\

\noindent \textbf{Remark. }A solution of this problem will imply the truth
of a certain number theoretic conjecture due to Levine and O'Sullivan \cite
{LOS} .

\section{Partial Solution}

For fixed $n$, let $\xi =$ $\left( x_1,x_2,...,x_n\right) $, 
\[
f(\xi )=\sum_{i=1}^n\frac 1{x_i} 
\]

\noindent and 
\[
P=\left\{ \xi :(j+1)x_j+x_i\geq (j+1)i+\varepsilon _{i\,j},1\leq j\leq i\leq
n\right\} 
\]

\noindent If $\xi \in P$, then clearly $x_i\geq 1$ for $1\leq i\leq n$. As a
consequence, $P$ contains no lines and $f$ is bounded above on $P$;
therefore, the supremum of $f$ over $P$ is attained at one or more vertices
of $P$ \cite{Rockafellar}. This proves the existence part of (i). While we
do not know how to prove the remainder of (i) or (ii), we show here that the
conjectured global maximum $\alpha =\left( a_1,a_2,...,a_n\right) $ is:

\begin{enumerate}
\item[(a)]  well-defined

\item[(b)]  feasible (that is, $\alpha \in P$)

\item[(c)]  a vertex of $P$
\end{enumerate}

\noindent and, if certain key inequalities hold,

\begin{enumerate}
\item[(d)]  a \textit{local} maximum of $f(\xi )$ subject to $\xi \in P$.
\end{enumerate}

\subsection{Proof of (a)}

The well-definition issue arises because of the conceivable non-uniqueness
(or even non-existence) of $j$ when determining $a_i$ for $i\geq 4$. Let 
\[
c_j=(j+2)a_{j+1}-(j+1)a_j 
\]

\noindent for $j\geq 1$. Expressed using $c_1,c_2,...$, the definition of $%
a_4,a_5...$ is 
\[
a_i=\left\{ 
\begin{tabular}{llll}
$3(i-2)$ &  & for & $4=c_1\leq i<c_2=10$ \\ 
$4(i-4)$ &  & for & $10=c_2\leq i<c_3=14$ \\ 
$5(i-6)$ &  & for & $14=c_3\leq i<c_4=24$ \\ 
$.$ &  &  &  \\ 
$.$ &  &  & 
\end{tabular}
\right. 
\]

\noindent We prove that both sequences $a_1,a_2,...$ and $c_1,c_2,...$ are
strictly increasing. Hypothesize inductively that $%
c_{i-1}>c_{i-2}>...>c_2>c_1=4$, where $i\geq 4$ is fixed. Since $c_{i-1}>i$,
there exists uniquely $j<i$ with $c_{j-1}\leq i<c_j$. If $i+1<c_j$, then $%
c_{j-1}\leq i+1<c_j$ and hence 
\[
a_{i+1}-a_i=(j+1)(i+1-a_j)-(j+1)(i-a_j)=j+1>0 
\]

\noindent If $i+1=c_j$, then $c_j\leq i+1<c_{j+1}$ and hence 
\begin{eqnarray*}
a_{i+1}-a_i &=&(j+2)(i+1-a_{j+1})-(j+1)(i-a_j) \\
&=&\left[ (j+2)(i+1)-(j+1)\,i\right] -\left[ (j+2)a_{j+1}-(j+1)a_j\right] \\
&=&\left[ (j+2)(i+1)-(j+1)\,i\right] -c_j=\left[ (j+2)(i+1)-(j+1)\,i\right]
-(i+1)=j+1>0
\end{eqnarray*}

\noindent We deduce that $(a_{i+1}-a_i)-(a_i-a_{i-1})\geq 0$ and thus 
\begin{eqnarray*}
c_i-c_{i-1}
&=&(i+2)a_{i+1}-2(i+1)a_i+i\,a_{i-1}=2a_{i+1}+i(a_{i+1}+a_{i-1})-2(i+1)a_i \\
&\geq &2a_{i+1}+2\,i\,a_i-2(i+1)a_i=2(a_{i+1}-a_i)>0
\end{eqnarray*}
\noindent This completes the inductive proof, from which well-definition
follows immediately. As a consequence, we may define 
\[
b_i=\left\{ 
\begin{tabular}{llll}
$1$ &  & if & $1\leq i\leq 3$ \\ 
$j$ &  & if & $4\leq i\leq n$ and $c_{j-1}\leq i<c_j$%
\end{tabular}
\right. 
\]
\noindent without ambiguity.

\subsection{Proof of (b)}

This is trivial if $1\leq i\leq 3$. If $i\geq 4$ and $j=b_i$, it follows
that 
\[
(j+1)(i-a_j)-(j+k+1)(i-a_{j+k})=\sum_{m=0}^{k-1}(c_{j+m}-i)>0 
\]
\noindent for $1\leq k\leq n-j$ and 
\[
(j+1)(i-a_j)-(j-k+1)(i-a_{j-k})=\sum_{m=1}^k(i-c_{j-m})\geq 0 
\]
\noindent for $0\leq k\leq j-1$. Both series are telescoping and the
inequalities are consequences of part (a). We deduce that 
\[
a_i=(j+1)(i-a_j)\geq (p+1)(i-a_p) 
\]
\noindent for any $1\leq p\leq n$, from which feasibility of $\alpha $
follows immediately.

\subsection{Proof of (c)}

This is true since $\alpha $ lies at the intersection of the $n$ hyperplanes 
\[
\begin{array}{ccc}
(b_i+1)x_{b_i}+x_i=(b_i+1)i+\varepsilon _{i\,b_i} &  & 1\leq i\leq n
\end{array}
\]

\subsection{Key Inequalities}

Before discussing part (d), we need to state certain key inequalities which,
although unproven, appear to be true for all $n\leq 10000$ via computer
check.\\

\noindent \textbf{Definition.} Fix integers $i$ and $j$ with $1\leq i<j$.
Let 
\begin{eqnarray*}
k_1 &=&j-1 \\
k_2 &=&b_{k_1} \\
k_3 &=&b_{k_2} \\
&&\ . \\
&&\ . \\
k_{m-1} &=&b_{k_{m-2}} \\
k_m &=&b_{k_{m-1}}
\end{eqnarray*}
\noindent where $m$ is the smallest integer such that $i\geq k_m$. Clearly
such an integer $m$ exists. Then define 
\[
d_{i\,j}=\left\{ 
\begin{tabular}{llll}
$(-1)^m\dprod\limits_{p=1}^m(k_p+1)$ &  & if & $i=k_m$ \\ 
$0$ &  & if & $i>k_m$%
\end{tabular}
\right. 
\]
\noindent \textbf{Conjecture.} Let $c_0=2$ for convenience, then 
\[
\begin{array}{cccc}
x_i^{*}\equiv \dfrac 1{a_i^2}+\dsum\limits_{j=i+1}^{b_n+1}\left(
d_{i\,j}\dsum\limits_{k=c_{j-2}}^{\min \{c_{j-1}-1,n\}}\dfrac
1{a_k^2}\right) \geq 0 &  & \text{for} & 1\leq i\leq n
\end{array}
\]
\noindent For example, if $n=24$, then $\xi ^{*}$ is the solution of the
linear system $M\,\xi =v$, where $M$ is the $24\times 24$ identity matrix
plus some upper triangular entries in the first $b_n=5$ rows as indicated: 
\[
\left( 
\begin{tabular}{llllllllllllllllllllllll}
{\tiny 1} & {\tiny 2} & {\tiny 2} & {\tiny 0} & {\tiny 0} & {\tiny 0} & 
{\tiny 0} & {\tiny 0} & {\tiny 0} & {\tiny 0} & {\tiny 0} & {\tiny 0} & 
{\tiny 0} & {\tiny 0} & {\tiny 0} & {\tiny 0} & {\tiny 0} & {\tiny 0} & 
{\tiny 0} & {\tiny 0} & {\tiny 0} & {\tiny 0} & {\tiny 0} & {\tiny 0} \\ 
& {\tiny 1} & {\tiny 0} & {\tiny 3} & {\tiny 3} & {\tiny 3} & {\tiny 3} & 
{\tiny 3} & {\tiny 3} & {\tiny 0} & {\tiny 0} & {\tiny 0} & {\tiny 0} & 
{\tiny 0} & {\tiny 0} & {\tiny 0} & {\tiny 0} & {\tiny 0} & {\tiny 0} & 
{\tiny 0} & {\tiny 0} & {\tiny 0} & {\tiny 0} & {\tiny 0} \\ 
&  & {\tiny 1} & {\tiny 0} & {\tiny 0} & {\tiny 0} & {\tiny 0} & {\tiny 0} & 
{\tiny 0} & {\tiny 4} & {\tiny 4} & {\tiny 4} & {\tiny 4} & {\tiny 0} & 
{\tiny 0} & {\tiny 0} & {\tiny 0} & {\tiny 0} & {\tiny 0} & {\tiny 0} & 
{\tiny 0} & {\tiny 0} & {\tiny 0} & {\tiny 0} \\ 
&  &  & {\tiny 1} & {\tiny 0} & {\tiny 0} & {\tiny 0} & {\tiny 0} & {\tiny 0}
& {\tiny 0} & {\tiny 0} & {\tiny 0} & {\tiny 0} & {\tiny 5} & {\tiny 5} & 
{\tiny 5} & {\tiny 5} & {\tiny 5} & {\tiny 5} & {\tiny 5} & {\tiny 5} & 
{\tiny 5} & {\tiny 5} & {\tiny 0} \\ 
&  &  &  & {\tiny 1} & {\tiny 0} & {\tiny 0} & {\tiny 0} & {\tiny 0} & 
{\tiny 0} & {\tiny 0} & {\tiny 0} & {\tiny 0} & {\tiny 0} & {\tiny 0} & 
{\tiny 0} & {\tiny 0} & {\tiny 0} & {\tiny 0} & {\tiny 0} & {\tiny 0} & 
{\tiny 0} & {\tiny 0} & {\tiny 6} \\ 
&  &  &  &  & {\tiny 1} &  &  &  &  &  &  &  &  &  &  &  &  &  &  &  &  &  & 
\\ 
&  &  &  &  &  & {\tiny 1} &  &  &  &  &  &  &  &  &  &  &  &  &  &  &  &  & 
\\ 
&  &  &  &  &  &  & {\tiny .} &  &  &  &  &  &  &  &  &  &  &  &  &  &  &  & 
\\ 
&  &  &  &  &  &  &  & {\tiny .} &  &  &  &  &  &  &  &  &  &  &  &  &  &  & 
\end{tabular}
\right) 
\]

\noindent and $v$ is the 24-vector with $i^{\text{th}}$ element $1/a_i^2$.
The inverse, $M^{-1}$, of $M$ is given by 
\[
\left( 
\begin{tabular}{llllllllllllllllllllllll}
{\tiny 1} & {\tiny -2} & {\tiny -2} & {\tiny 6} & {\tiny 6} & {\tiny 6} & 
{\tiny 6} & {\tiny 6} & {\tiny 6} & {\tiny 8} & {\tiny 8} & {\tiny 8} & 
{\tiny 8} & {\tiny -30} & {\tiny -30} & {\tiny -30} & {\tiny -30} & {\tiny %
-30} & {\tiny -30} & {\tiny -30} & {\tiny -30} & {\tiny -30} & {\tiny -30} & 
{\tiny -36} \\ 
& {\tiny 1} & {\tiny 0} & {\tiny -3} & {\tiny -3} & {\tiny -3} & {\tiny -3}
& {\tiny -3} & {\tiny -3} & {\tiny 0} & {\tiny 0} & {\tiny 0} & {\tiny 0} & 
{\tiny 15} & {\tiny 15} & {\tiny 15} & {\tiny 15} & {\tiny 15} & {\tiny 15}
& {\tiny 15} & {\tiny 15} & {\tiny 15} & {\tiny 15} & {\tiny 18} \\ 
&  & {\tiny 1} & {\tiny 0} & {\tiny 0} & {\tiny 0} & {\tiny 0} & {\tiny 0} & 
{\tiny 0} & {\tiny -4} & {\tiny -4} & {\tiny -4} & {\tiny -4} & {\tiny 0} & 
{\tiny 0} & {\tiny 0} & {\tiny 0} & {\tiny 0} & {\tiny 0} & {\tiny 0} & 
{\tiny 0} & {\tiny 0} & {\tiny 0} & {\tiny 0} \\ 
&  &  & {\tiny 1} & {\tiny 0} & {\tiny 0} & {\tiny 0} & {\tiny 0} & {\tiny 0}
& {\tiny 0} & {\tiny 0} & {\tiny 0} & {\tiny 0} & {\tiny -5} & {\tiny -5} & 
{\tiny -5} & {\tiny -5} & {\tiny -5} & {\tiny -5} & {\tiny -5} & {\tiny -5}
& {\tiny -5} & {\tiny -5} & {\tiny 0} \\ 
&  &  &  & {\tiny 1} & {\tiny 0} & {\tiny 0} & {\tiny 0} & {\tiny 0} & 
{\tiny 0} & {\tiny 0} & {\tiny 0} & {\tiny 0} & {\tiny 0} & {\tiny 0} & 
{\tiny 0} & {\tiny 0} & {\tiny 0} & {\tiny 0} & {\tiny 0} & {\tiny 0} & 
{\tiny 0} & {\tiny 0} & {\tiny -6} \\ 
&  &  &  &  & {\tiny 1} &  &  &  &  &  &  &  &  &  &  &  &  &  &  &  &  &  & 
\\ 
&  &  &  &  &  & {\tiny 1} &  &  &  &  &  &  &  &  &  &  &  &  &  &  &  &  & 
\\ 
&  &  &  &  &  &  & {\tiny .} &  &  &  &  &  &  &  &  &  &  &  &  &  &  &  & 
\\ 
&  &  &  &  &  &  &  & {\tiny .} &  &  &  &  &  &  &  &  &  &  &  &  &  &  & 
\end{tabular}
\right) 
\]

\noindent and the entries $x_i^{*}$ of $\xi ^{*}=M^{-1}v$ are prescribed by
the above summation formula. In the case $n=24$, we compute 
\[
x_1^{*}=\frac{123587941503427}{187646731272000}>0 
\]
\[
x_2^{*}=\frac{3536905093973}{27799515744000}>0 
\]
\[
x_3^{*}=\frac{44159}{1016064}>0 
\]

\[
x_4^{*}=\frac{9439261073843}{750586925088000}>0 
\]
\[
x_5^{*}=\frac{47}{4050}>0 
\]

\noindent and these positivity results are consistent with the Conjecture.
Of course, $x_i^{*}>0$ for $i>b_n$ immediately.

\subsection{Partial proof of (d)}

It suffices to solve the following (\textit{primal}) linear programming
problem: 
\[
\begin{array}{ccc}
\text{Minimize} &  & g(\xi )=(-\xi )\cdot \nabla f(\alpha
)=\dsum\limits_{i=1}^n\dfrac{x_i}{a_i^2}
\end{array}
\]
\[
\begin{array}{ccc}
\text{subject to} &  & \xi \text{ }\epsilon \text{ }Q
\end{array}
\]

\noindent where $Q$ is the polyhedron 
\[
Q=\left\{ \xi :(b_i+1)x_{b_i}+x_i\geq (b_i+1)i+\varepsilon _{i\,b_i},1\leq
i\leq n\right\} 
\]

\noindent Note that $Q$ contains $P$ and possesses a unique vertex, $\alpha $%
. Note also that, by the Conjecture, the \textit{dual} linear programming
problem has nonempty feasible region 
\[
R=\left\{ 
\begin{array}{cccccc}
\xi : &  & x_j+(j+1)\dsum\limits_{i=c_{j-1}}^{\min \{c_j-1,n\}}x_i=\dfrac
1{a_j^2} &  & \text{for} & 1\leq j\leq b_n \\ 
&  & x_j=\dfrac 1{a_j^2} &  & \text{for} & b_n<j\leq n \\ 
&  & x_j\geq 0 &  & \text{for} & 
\begin{tabular}{ll}
all & $j$%
\end{tabular}
\end{array}
\right\} =\left\{ M^{-1}v\right\} =\left\{ \xi ^{*}\right\} 
\]

\noindent hence $g$ is bounded below on $Q$. Therefore $\alpha $ is the
global minimum of $g(\xi )$ subject to $\xi $ $\epsilon $ $Q$, which implies
that $\alpha $ is a local maximum of $f(\xi )$ subject to $\xi $ $\epsilon $ 
$P$.

\subsection{Partial proof of the Conjecture}

The key inequalities are provably true when $i$ is sufficiently large
relative to $n$. More precisely, if 
\[
b_{b_{b_n}}<i\leq n 
\]
\noindent then 
\[
x_i^{*}\geq \dfrac 1{a_i^2}-(i+1)\dsum\limits_{j=c_{i-1}}^{c_i-1}\dfrac
1{a_j^2}>0 
\]
\noindent To see this, we prove two lemmas.\\

\noindent \textbf{Lemma One. } $d_{i,\,i+1}=-(i+1)$ for all $i\geq 1$ and $%
d_{i\,j}\geq 0$ if $i>b_{b_{b_n}}$ and $i+1<j\leq b_n+1$.\\

\noindent \textbf{Proof of Lemma One. }The first part is trivial. The second
part is proved by noting that $m>1$ since $k_1=j-1>i$, so either $m=2$
(which implies that $d_{i\,j}\geq 0$) or $m=3$ since 
\[
k_3=b_{b_{j-1}}\leq b_{b_{b_n}}<i 
\]

\noindent  (which, in turn, implies that $d_{i\,j}=0$). QED.\\

\noindent \textbf{Lemma Two. } $\dfrac
1{a_j^2}-(j+1)\dsum\limits_{i=c_{j-1}}^{c_j-1}\dfrac 1{a_i^2}>0$ for all $%
j\geq 1$.\\

\noindent \textbf{Proof of Lemma Two. }Direct computation proves the
inequality for $j=1,2,3,4,10$ and $14$. For all other values of $j$, we will
show that 
\[
\dfrac 1{a_j^2}-\dfrac{c_j-c_{j-1}}{j+1}\dfrac 1{(c_{j-1}-a_j)^2}>0 
\]
\noindent that is, 
\[
e_j\equiv (j+1)c_{j-1}(c_{j-1}-2\,a_j)+\left[ (j+1)-(c_j-c_{j-1})\right]
a_j^2>0 
\]
\noindent which implies the truth of the Lemma. Observe that, if $%
c_{k-1}<j\leq c_k$, then 
\[
c_{j-1}-2\,a_j=(k+1)a_k>0 
\]
\noindent and 
\[
(j+1)-(c_j-c_{j-1})=\left\{ 
\begin{tabular}{llll}
$j-2\,k-1\geq 0$ &  & if & $c_{k-1}<j<c_k$ \\ 
$-2\,k-3<0$ &  & if & $j=c_k$%
\end{tabular}
\right. 
\]
\noindent These inequalities yield $e_j>0$ when $j\neq c_k$ for any k. In
the event $j=c_k$ for some $k\geq 4$, the argument is only slightly more
complicated: 
\begin{eqnarray*}
e_j &=&(k+1)^2\left[ (c_k+1)(2\,c_k-a_k)a_k-(2\,k+3)(c_k-a_k)^2\right] \\
&\geq &(k+1)^2\left[ (c_k-a_k)(2\,c_k-2\,a_k)a_k-(2\,k+3)(c_k-a_k)^2\right]
\\
&=&(k+1)^2(c_k-a_k)^2\left[ 2\,a_k-(2\,k+3)\right] >0
\end{eqnarray*}
\noindent for all $k\geq 4$. QED.

\section{Closing Words}

Techniques for numerical convex maximization abound \cite{PR}. A vertex
enumeration scheme has led to verification that $\alpha $ is the \textit{%
global} maximum of $f(\xi )$ subject to $\xi $ $\epsilon $ $P$ for small $n$
only. Keith Briggs has used the general-purpose optimization programs 
\textsc{AMPL} and \textsc{LANCELOT} to confirm the global maximum claim up
to $n=24$, and \textsc{CFSQP} to do likewise up to $n=121$.\\

The continuous analog of this problem (with summations replaced by
integrals) is discussed in a companion paper.\\

An outcome of Levine and O'Sullivan's work \cite{LOS} is that, for any $n$,
there is a global maximum $\alpha $ that satisfies $a_1=1,$ $a_2=2,$ $a_3=4$
and either $a_4=6$ or $a_4\geq 28,$ where $\xi $ is restricted to \textit{%
integer} points in $P$ (that is, to $\xi $ $\epsilon $ $P\cap \Bbb{Z}^n$).
Their proof unfortunately does not extend to the real case.\\

Do there exist other functions $f$ and polyhedra $P$ for which the
maximizing vertex $\alpha $ is ''self-generating'' as the dimension $n$
increases? A simple characterization of such pairs $(f,P)$ may lead to the
insight necessary to solve this problem.\newpage\

\end{document}